\documentstyle[11pt]{article}

\makeatletter

\@addtoreset{equation}{section}
\makeatother
\def\<{\langle}
\def\>{\rangle}

\def\epp{\varepsilon}

\newtheorem{lem}{Lemma}[section]
\newtheorem{theo}{Theorem}[section]
\newtheorem{rem}{Remark}[section]

\makeatletter
   
   \@addtoreset{equation}{section}
\makeatother

\begin{document}
\title{\bf Asymptotic profiles for damped plate equations \\ with rotational inertia terms}
\author{Tomonori Fukushima and Ryo Ikehata\thanks{Corresponding author: ikehatar@hiroshima-u.ac.jp} \\{\small Department of Mathematics} \\{\small Graduate School of Education} \\ {\small Hiroshima University} \\ {\small Higashi-Hiroshima 739-8524, Japan} \\and\\Hironori Michihisa \\ {\small Department of Mathematics}\\ {\small Graduate School of Science} \\ {\small Hiroshima University} \\ {\small Higashi-Hiroshima 739-8526, Japan}}
\maketitle
\begin{abstract}
We consider the Cauchy problem for plate equations with rotational inertia and frictional damping terms. We will derive asymptotic profiles of the solution in $L^{2}$-sense as $t \to \infty$ in the case when the initial data have high and low regularity, respectively. Especially, in the low regularity case of the initial data one encounters the regularity-loss structure of the solutions, and the analysis is more delicate. We employ the so-called Fourier splitting method combined with the explicit expression of the solution (high frequency estimates) and the method due to \cite{Ike-1} (low frequency estimates).  
\end{abstract}
\section{Introduction}
\footnote[0]{Keywords and Phrases: Plate equation; frictional damping; rotational inertia terms; asymptotic profiles, low frequency, high frequency.}
\footnote[0]{2010 Mathematics Subject Classification. Primary 35L05; Secondary 35B40, 35L30.}
We are concerned with the Cauchy problem of the following double dispersion equation with rotational inertia and frictional damping terms: 
\begin{equation}
\left\{\begin{array}{ll}
u_{tt}(t,x)-\Delta u_{tt}(t,x) + \Delta^{2}u(t,x) -\Delta u(t,x) + u_{t}(t,x) = 0, \quad (t,x) \in (0,\infty)\times  {\bf R}^{n}, \\[2mm]
\hspace{1cm} u(0,x)= u_{0}(x), \quad u_{t}(0,x) = u_{1}(x), \quad x \in {\bf R}^{n},
\end{array}
\right.
\end{equation}
where $n\ge1$. Here, one assumes for the moment
\[[u_{0},u_{1}] \in H^{3}({\bf R}^{n})\times H^{2}({\bf R}^{n}).\]
Under these conditions it is already known that problem (1.1) has a unique solution (see \cite[Theorem 2.1]{GHC}) in the class
\[u \in X_{2} := C^{2}([0,\infty); H^{1}({\bf R}^{n}))\cap C^{1}([0,\infty); H^{2}({\bf R}^{n}))\cap C([0,\infty); H^{3}({\bf R}^{n})).\]

The equation of (1.1) can be viewed as the wave equation with a weaker dissipative term:
\begin{equation}
u_{tt}(t,x)-\Delta u(t,x) + (1-\Delta)^{-1}u_{t}(t,x) = 0, 
\end{equation} 
which is proposed in the early papers \cite{SK, D, D-1}. In fact, \cite{SK, D, D-1} study a more generalized system of viscoelasticity of (1.2):   
\begin{equation}
u_{tt}-\sum_{j}b^{j}(\partial_{x}u)_{x_{j}} + (1-\Delta)^{-1}\left(\sum_{j,k}K^{jk}*u_{x_{j}x_{k}} + Lu_{t}\right) = 0.
\end{equation}
Here, the terms $\sum_{j,k}K^{jk}*u_{x_{j}x_{k}}$ and $Lu_{t}$ represent the dissipations of the memory-type and the friction type, respectively. If the memory effect vanishes, then the system (1.3) corresponds to the equation (1.2) as the special case of (1.3). In \cite{SK, D, D-1}, they study sharp decay estimates of the solutions to (1.3) based on the time-weighted energy method in order to overcome some difficulties coming from the so-called regularity-loss type property. 

While, quite recently, the following equation is studied in \cite{DL}:  
\begin{equation}
u_{tt}(t,x)-\Delta u(t,x) + (1-\Delta)^{-1}(-\Delta u_{t})(t,x) = 0, 
\end{equation} 
which includes the structural damping in place of the friction term. Equation (1.4) can be written as
\begin{equation}
u_{tt}(t,x)-\Delta u_{tt}(t,x) + \Delta^{2}u(t,x) -\Delta u(t,x) -\Delta u_{t}(t,x) = 0.
\end{equation}
In \cite{DL} they have studied long time decay estimates for the solution in $L^{p}$ spaces and in real Hardy spaces ${\cal H}^{1}$, and have applied them to the nonlinear problem with the nonlinearity such as $\Delta f(u)$. However, it seems that in \cite{DL} they still study neither the asymptotic profile of the solution to (1.4) nor the frictional damping case (1.1) itself, let alone the asymptotic profile of the solution to problem (1.1). 

In this connection, the asymptotic profile and the optimal decay of the solution to the generalized equation
\begin{equation}
u_{tt}(t,x) + (-\Delta)^{\delta}u_{tt}(t,x) + b\Delta^{2}u(t,x) +a(-\Delta)^{\alpha}u(t,x) + (-\Delta)^{\theta}u_{t}(t,x) = 0
\end{equation}
have been already investigated in \cite{HIC} in the case when $0 < \delta < \theta$ and $\frac{1}{2} < \theta \min\{\frac{3}{2}, \delta + \frac{1}{2}\}$. For the case of $\delta = 1 = b$, $a = 0$, and $\theta \in [0,1]$, \cite{CLI} studies almost optimal decay order of the total energy by developing a new type of energy method in the Fourier space. One can also cite \cite{TY} for $\delta = 0$, $\alpha = 1$ and $\theta = 0$ case, and it is also interesting to cite the paper \cite{DCL}, which studies (almost) optimal total energy decay rates of the equations (1.6) with $\delta = b =1$, $a =0$ and the fractional term $(-\Delta)^{\theta}u_{t}$ replaced by the time dependent fractional damping $2b(t)(-\Delta)^{\theta}u_{t}$. Furthermore, \cite{LC} with $\delta = b = \theta = 1$ and $a = 0$ in (1.6) derives various decay estimates of solutions together with application to nonlinear problems in exterior domains. Anyway, unfortunately the results in \cite{HIC} do not include the case (1.1) concerning the study of asymptotic profiles and optimality of decay rates of the solution. Problem (1.1) corresponds to the case for $\delta = \alpha = 1$ and $\theta = 0$ in (1.6). The reason for this difficulty comes from the fact that in order to study the equation (1.1) one has to get the more essential information in the high frequency region of the Fourier transformed solutions, which encounter the so-called regularity-loss structure. This regularity-loss property is extremely difficult when one catches the leading term of asymptotic expansions of the solution. In this sense, at least in the equation (1.1) case this is the first trial to investigate the asymptotic profile of solutions in the regularity-loss case.     

The main purpose of this paper is to catch the asymptotic profiles of the solution as $t \to \infty$ to problem (1.1), which is not well studied yet, and is to investigate the optimal decay rates of $L^2$-norm of solutions. 

This paper is organized as follows. In section 2, one prepares several notation, including the explicit formula of the Fourier transform to the solution. We state main results in Section 3. Section 4 will be devoted to proofs of our main results. Optimal decay estimates of solutions in $L^{2}$-norm will be given in section 5.


\section{Notation.}

Throughout this paper, we handle with the following function spaces and norms:
\[ 
\|f\|_{p}:=\left(\int_{{\bf R}^{n}} |f(x)|^p dx \right)^{\frac{1}{p}},
\hspace{5mm} 
f \in L^p({\bf R}^{n}); 
\]
\[ 
f \in L^{1,1}({\bf R}^{n}) 
\Leftrightarrow 
f \in L^1({\bf R}^{n}), \ 
\|f\|_{1,1}:=\int_{{\bf R}^{n}} (1+|x|)\vert f(x)\vert dx<\infty. 
\]
We also denote the Fourier transform $\hat{f}$ of the function $f$\ by
\[
\hat{f}(\xi)
:=\int_{{\bf R}^{n}} f(x)\ e^{-i\xi \cdot x} dx,
\]
where $i := \sqrt{-1}$. As usual, the norm of the Sobolev spaces $H^l ({\bf R}^{n})$ ($l \geq 0$) is written by 
\[ 
\|f\|_{H^{l}}
:=\left(\int_{{\bf R}^{n}} (1+|\xi|^2)^l |\hat{f}(\xi)|^2 d\xi \right)^{\frac{1}{2}}, 
\hspace{5mm} 
f \in H^l({\bf R}^{n}). 
\]
Taking the Fourier transform of the both sides of (1.1), one has 
\begin{equation}
\left\{\begin{array}{l}
(1+|\xi|^2)\hat{u}_{tt}+\hat{u}_{t}+|\xi|^2(1+|\xi|^2) \hat{u}=0 \,,  \\[2mm]
\hat{u}(0,\xi)=\hat{u}_{0},
\quad
\hat{u}_{t}(0,\xi)=\hat{u}_{1}(\xi).
\end{array}
\right.
\end{equation}
The characteristic equation of (2.1) is 
\[
(1+|\xi|^2)\lambda^2+\lambda+|\xi|^2(1+|\xi|^2)=0,
\]
and we set the solutions as 
\[ \lambda_{1}:=\frac{-1+\sqrt{1-4|\xi|^2(1+|\xi|^2)^2}}{2(1+|\xi|^2)} \,, \hspace{5mm} \lambda_{2}:=\frac{-1-\sqrt{1-4|\xi|^2(1+|\xi|^2)^2}}{2(1+|\xi|^2)} \,. \]
Thus the solution to (2.1) is given by
\begin{eqnarray}
\hat{u}(t,\xi)&=&\frac{e^{\lambda_{1}t}-e^{\lambda_{2}t}}{\lambda_{1}-\lambda_{2}}\hat{u}_{1}(\xi)+\frac{\lambda_{1}e^{\lambda_{2}t}-\lambda_{2}e^{\lambda_{1}t}}{\lambda_{1}-\lambda_{2}}\hat{u}_{0}(\xi) \\[2mm]
&=&\hat{u}_{0}(\xi) \cdot E_{0}(t,\xi)+\left\{\hat{u}_{1}(\xi)+\frac{1}{2(1+|\xi|^2)} \cdot \hat{u}_{0}(\xi) \right\} E_{1}(t,\xi),
\end{eqnarray}
where
\[
E_{0}(t,\xi):=\left\{
\begin{array}{ll}
e^{-\frac{t}{2(1+|\xi|^2)}} \cdot \displaystyle{\cos{\left(\frac{t\sqrt{4|\xi|^2(1+|\xi|^2)^2-1}}{2(1+|\xi|^2)} \right)}} &(|\xi| > \zeta), \\[6mm]
e^{-\frac{t}{2(1+|\xi|^2)}} \cdot \displaystyle{\cosh{\left(\frac{t\sqrt{1-4|\xi|^2(1+|\xi|^2)^2}}{2(1+|\xi|^2)} \right)}} &(|\xi|\le\zeta),
\end{array}
\right.
\]
\[
E_{1}(t,\xi):=\left\{
\begin{array}{ll}
e^{-\frac{t}{2(1+|\xi|^2)}} 
\left[\displaystyle{\sin{\left(\frac{t\sqrt{4|\xi|^2(1+|\xi|^2)^2-1}}{2(1+|\xi|^2)} \right)} \biggr/ \frac{\sqrt{4|\xi|^2(1+|\xi|^2)^2-1}}{2(1+|\xi|^2)}} \right] &(|\xi| > \zeta), \\[6mm]
e^{-\frac{t}{2(1+|\xi|^2)}} 
\left[\displaystyle{\sinh{\left(\frac{t\sqrt{1-4|\xi|^2(1+|\xi|^2)^2}}{2(1+|\xi|^2)} \right)} \biggr/ \frac{\sqrt{1-4|\xi|^2(1+|\xi|^2)^2}}{2(1+|\xi|^2)}} \right] &(|\xi|\le\zeta).
\end{array}
\right.
\]
Here, $\zeta \in (0,1)$ is a constant satisfying
\begin{equation}
4\zeta^2(1+\zeta^2)^2-1=0.
\end{equation}
For later use, one defines
\[P_{j} := \int_{{\bf R}^{n}}u_{j}(x)dx\quad (j = 0,1).\]


\section{Main results}
In this section, we state our main results of this paper. The leading terms of the solution shown in Theorems~\ref{thm:1} and \ref{thm:2} below are obtained by expanding the solution formula (2.3) in the high frequency region and the low frequency region, respectively. 
These results indicate that the solution behaves like the Gauss kernel under the high regularity condition, while it has strong wave-like property for low regularity initial data.

\begin{theo}
\label{thm:1}
Let $n \ge 7$. 
If $(u_{0},u_{1}) \in (H^{l+1}({\bf R}^{n}) \cap L^{1,1}({\bf R}^{n})) \times (H^{l}({\bf R}^{n}) \cap L^{1,1}({\bf R}^{n}))$, 
then the solution $u \in X_{2}$ to {\rm (1.1)} satisfies 
\begin{eqnarray}
& &\left\|\hat{u}(t,\xi)-\left(\hat{u}_{1}(\xi) \cdot e^{-\frac{t}{2|\xi|^2}} \frac{\sin{(t|\xi|)}}{|\xi|}+\hat{u}_{0}(\xi) \cdot e^{-\frac{t}{2|\xi|^2}}\cos{(t|\xi|)} \right) \right\|_2  \nonumber \\[2mm]
&\le& \left\{ \begin{array}{ll}
CI_0\, t^{-\frac{l+3}{2}}, 
& \,\, 
2\le l\leq n/2-3 
\,\,\,
\mbox{with}
\,\,\,
n\ge10, \\[2mm]
CI_0\, t^{-\frac{n}{4}},
&\,\, 
2\le l <n/2 -1 
\,\,\,
\mbox{with}
\,\,\,
n=7, 8, 9\\[2mm] & \hspace{3cm}\mbox{or}\,\,\, 
n/2-3 <l < n/2-1
\,\,\,\mbox{with}\,\,\,n\ge10, \\[0.2cm]
\end{array}
\right. \nonumber
\end{eqnarray}
for $t \ge 1$. 
Here, we put 
\[
I_0:=\|u_1\|_{H^l}+\|u_0\|_{H^{l+1}} + \| u_{1} \|_{1,1} + \| u_{0}\|_{1,1}
\]
and $C>0$ is a constant independent of $t$ and $u_0,\,u_1$.
\end{theo}

\begin{rem}{\rm (i) The diffusion wave property like
\begin{equation}
\hat{u}(t,\xi) \sim \hat{u}_{1}(\xi) \cdot e^{-\frac{t\vert\xi\vert^{2}}{2}} \frac{\sin{(t|\xi|)}}{|\xi|}+\hat{u}_{0}(\xi) \cdot e^{-\frac{t\vert\xi\vert^{2}}{2}}\cos{(t|\xi|)}
\end{equation}
in asymptotic sense has been already discussed in previous papers\\
by \cite{Ike-1, ITY} to the equation
\[u_{tt}-\Delta u -\Delta u_{t} = 0,\]
by \cite{IS} to the equation 
\[u_{tt} + \Delta^{2} u -\Delta u -\Delta u_{t} = 0,\]
and by \cite{WC} to the equation
\[u_{tt} -\Delta u_{tt} + \Delta^{2} u -\mu^{2}\Delta u -\nu\Delta u_{t} = 0,\]
respectively. It should be noticed that the profile in asymptotic sense derived in Theorem 3.1 is different from that of (3.1) discovered in \cite{Ike-1, ITY, IS, WC}. In this connection, Michihisa \cite{Michi} first derived such type of wave-like profiles stated in Theorem 3.1 in the case of low regularity solutions toward the so-called Rosenau equation:
\[u_{tt} - \Delta u - \nu\Delta u_{t} + \Delta^{2} u + \Delta^{2} u_{tt} = 0.\]
The Rosenau equation above also includes a regularity-loss structure in the high frequency region of the solutions. \\
(ii) In the low regularity case $2 \le l\le n/2-3$, analysis in the high frequency region is essential due to the regularity-loss effect. 
Even when we impose a little stronger regularity assumption such as $n/2-3< l < n/2-1$, 
we can expect the wave part 
\begin{equation}
\hat{u}_{1}(\xi) \cdot e^{-\frac{t}{2|\xi|^2}} \frac{\sin{(t|\xi|)}}{|\xi|}+\hat{u}_{0}(\xi) \cdot e^{-\frac{t}{2|\xi|^2}}\cos{(t|\xi|)}
\end{equation}
is still dominant in the asymptotic sense. 
Confirming function~(3.2) decays like $t^{-\frac{l+1}{2}}$ in the $L^2$ space, 
the constraint $l \le n/2-1$ assures that profile~(3.2) decays slower than the heat kernel. 
On the other hand, stronger smoothness $n/2-1<l$ makes it possible to simplify analysis in the high frequency region, for the diffusive structure is mainly observed (see Theorem~3.2 below). 
Related to this concern, see also Theorems~5.1 and 5.2.
}
\end{rem}

\begin{theo}
\label{thm:2}
Let $n \ge 1$. 
If $(u_{0},u_{1}) \in (H^{l+1}({\bf R}^{n}) \cap L^{1,1}({\bf R}^{n})) \times (H^{l}({\bf R}^{n}) \cap L^{1,1}({\bf R}^{n}))$, 
 then the solution $u \in X_{2}$ to {\rm (1.1)} satisfies 
\[\left\|\hat{u}(t,\xi)-(P_1+P_0)\, e^{-t|\xi|^2} \right\|_2\]
\[\leq \left\{
\begin{array}{ll}
C I_0\, t^{-\frac{l+1}{2}}, 
&  \quad  
\ell=2 \,\,\,\mbox{with}\,\,\,n=4
\quad\mbox{or}\quad
2\le l\le 5/2 \,\,\,\mbox{with}\,\,\,n=5\\[2mm]
& \hspace{3.5cm}\mbox{or}\,\,\, 
n/2-1 <l \leq n/2
\,\,\,\mbox{with}\,\,\,n\ge6, \\[0.2cm]
C I_0 \, t^{-\frac{n+2}{4}}, 
& \quad
n/2<l\,\,\,\mbox{with}\,\,\,n\ge4\\[2mm]
& \hspace{2cm}\mbox{or}\,\,\, 
2 \leq l
\,\,\,\mbox{with}\,\,\,n = 1,2,3,\\[0.2cm]
\end{array} \right. \]
for $t \ge 1$. 
Here, $C>0$ is a constant independent of $t$ and $u_0,\,u_1$.
\end{theo}
\begin{rem}{\rm In the case of high regularity solutions, one can see from Theorem 3.2 that the asymptotic profile of solutions to problem (1.1) is strongly dominated from the equation:
\[u_{tt} -\Delta u + u_{t} = 0,\]
which is quite natural by considering the limit $\vert\xi\vert \to 0$ in the Fourier transformed equation of (1.2):
\[\hat{u}_{tt}+\vert\xi\vert^{2}\hat{u} + (1+|\xi|^2)^{-1}\hat{u}_{t}=0.\] 
}

\begin{theo}
\label{thm:3}
Let $n \ge 6$ and $l=n/2-1$. 
If $(u_{0},u_{1}) \in (H^{l+1}({\bf R}^{n}) \cap L^{1,1}({\bf R}^{n})) \times (H^{l}({\bf R}^{n}) \cap L^{1,1}({\bf R}^{n}))$, 
 then the solution $u \in X_{2}$ to {\rm (1.1)} satisfies 
\begin{eqnarray}
& & \left\|\hat{u}(t,\xi)-\left\{(P_1+P_0)\, e^{-t|\xi|^2}
+\hat{u}_{1}(\xi) \cdot e^{-\frac{t}{2|\xi|^2}} \frac{\sin{(t|\xi|)}}{|\xi|}+\hat{u}_{0}(\xi) \cdot e^{-\frac{t}{2|\xi|^2}}\cos{(t|\xi|)} \right\} \right\|_2 \nonumber \\[2mm]
&\le& C I_0\, t^{-\frac{n+2}{4}} \nonumber
\end{eqnarray}
for $t \ge 1$. 
Here, $C>0$ is a constant independent of $t$ and $u_0,\,u_1$.
\end{theo}
\end{rem}
\begin{rem}{\rm The condition $l \geq 2$ in Theorems 3.1-3.3 implies that the problem (at least) has a unique solution $u \in X_{2}$. In fact, one has to restrict the regularity condition $\ell$ to satisfy the relation with $n$ when one uses Theorems above. }
\end{rem}

\section{Proofs of main results}
We prove in this section our main results stated in section 3. For this purpose, one first takes a constant $0<\delta<\zeta$ satisfying 
\[
\frac{1}{2} \le 1-4|\xi|^2(1+|\xi|^2)^2 \le 1
\quad 
\mbox{if}
\quad
|\xi| \le \delta.
\]
We prepare the following lemmas.  
\begin{lem}
Let $n \ge 1$ and $l \ge 0$. 
If $(u_{0},u_{1}) \in H^{l+1}({\bf R}^{n}) \times H^{l}({\bf R}^{n})$,  
then the function $u$ defined in (2.3) satisfies
\begin{eqnarray}
& &\int_{|\xi| \ge 1} \left|\hat{u}(t,\xi)-\left(\hat{u}_{1}(\xi) \cdot e^{-\frac{t}{2|\xi|^2}} \frac{\sin{(t|\xi|)}}{|\xi|}+\hat{u}_{0}(\xi) \cdot e^{-\frac{t}{2|\xi|^2}} \cos{(t|\xi|)} \right) \right|^2 d\xi \nonumber \\[2mm]
&\le& C\left(||u_{1}||_{H^l}^2+||u_{0}||_{H^{l+1}}^2 \right)\, (1+t)^{-l-3} \nonumber
\end{eqnarray}
for $t \ge 0$. Here, $C>0$ is a constant independent of $t$ and $u_0$, $u_1$.
\end{lem}
\begin{lem}
Let $n \ge 1$. 
If $(u_{0},u_{1}) \in L^{1,1}({\bf R}^{n}) \times L^{1,1}({\bf R}^{n})$, 
then the function $u$ defined in (2.3) satisfies
\begin{eqnarray}
\int_{|\xi| \le \delta} \left|\hat{u}(t,\xi)-(P_{1}+P_{0})\, e^{-t|\xi|^2} \right|^2 d\xi 
&\le& C\left(\|u_{1}\|_{1,1}^2+\|u_{0}\|_{1,1}^2\right)\, (1+t)^{-\frac{n+2}{2}} \nonumber
\end{eqnarray}
fot $t \ge 0$. Here, $C>0$ is a constant independent of $t$ and $u_0$, $u_1$.
\end{lem}

We can ignore the behavior of the solution in the middle frequency region, for its $L^2$-norm is exponentially small. 
\begin{lem}
Let $n \ge 1$. 
If $(u_{0},u_{1}) \in L^{2}({\bf R}^{n}) \times L^{2}({\bf R}^{n})$,  
then the function $u$ defined in (2.3) satisfies
\[ \int_{\delta<|\xi|<1} |\hat{u}(t,\xi)|^2 d\xi \le C \left(\|u_{0}\|_2^2+\|u_{1}\|_2^2 \right)\, e^{-\eta t} \]
for $t \ge 0$. Here, $C>0$ and $\eta>0$ are constants independent of $t$ and $u_0$, $u_1$. \\
\end{lem}

In order to prove above lemma 4.1, 
we prepare the following lemma. 
\begin{lem}
Let $l \ge 0$. 
Then there exists a constant $C>0$ such that 
\[
\sup_{x \ge 1}{\left(\frac{e^{-\frac{t}{x}}}{x^{l}}\right)} 
\le C(1+t)^{-l}
\]
for $t \ge 0$. 
\end{lem}
{\it Proof.} 
One easily sees that 
\[
\sup_{x \ge 1}{\left(\frac{e^{-\frac{t}{x}}}{x^{l}}\right)} 
\le e \sup_{x \ge 1}{\left(\frac{e^{-\frac{1+t}{x}}}{x^{l}}\right)}
\]
for $t\ge0$. 
Now we change variables 
\[
s:=\sqrt{\frac{1+t}{x}}, 
\quad
\mbox{i.e.,} 
\quad 
x=\frac{1+t}{s^2}.
\]
Then, it follows that 
\[
\sup_{x \ge 1}{\left(\frac{e^{-\frac{1+t}{x}}}{x^{l}}\right)}
=\sup_{0 \le s \le \sqrt{1+t}}{\left(s^{2l} \, e^{-s^2} \right)} \, (1+t)^{-l}
\le \sup_{s \ge 0}{\left(s^{2l} \, e^{-s^2} \right)} \, (1+t)^{-l}
\le C(1+t)^{-l},
\]
which implies the desired estimate. 
\hfill
$\Box$ 
\par
\vspace{0.1cm}

When one proves Lemma 4.2, 
the following decomposition of the initial data is available. 
It is introduced by Ikehata \cite{Ike-0}: 
\[ \hat{u}_{j}(\xi) =A_{j}(\xi)-iB_{j}(\xi)+P_{j}\quad (j=0,1), \]
where
\[A_{j}=A_{j}(\xi) := \int_{{\bf R}^{n}} u_{j}(x)\,\big(\cos{(\xi \cdot x)}-1\big) dx,\quad B_{j}=B_{j}(\xi) := \int_{{\bf R}^{n}} u_{j}(x)\,\sin{(\xi \cdot x)}dx,\]
\[P_{j} =\int_{{\bf R}^{n}} u_{j}(x)  dx. \]
Here, we recall the following useful lemma.
\begin{lem}
{\rm \cite{Ike-0}}
Let $n \ge 1$. 
Then it holds for all $\xi \in {\bf R}^{n}$,
\[ |A_{j}(\xi)| \le L\|u_{j}\|_{1,1}|\xi|,\quad |B_{j}(\xi)| \le M\|u_{j}\|_{1,1}|\xi|,\]
where
\[ L:=\sup_{\theta \neq 0}{\frac{|1-\cos{\theta}|}{|\theta|}}<+\infty,\quad 
M:=\sup_{\theta \neq 0}{\frac{|\sin{\theta}|}{|\theta|}}<+\infty. \]
Here, $A_j$ and $B_j$ are functions defined above for $u_j \in L^{1,1}({\bf R}^{n})$ {\rm(}$j=0,1${\rm)}.
\end{lem}

Let us prove Lemma 4.1 first.

{\it Proof of Lemma 4.1.}\,For simplicity, put $\alpha=\alpha(\xi) := 1+|\xi|^2$. 
By the mean value theorem one has 
\begin{eqnarray}
\sin{\left(\frac{t\sqrt{4\alpha^2 |\xi|^2-1}}{2\alpha}\right)}-\sin{(t|\xi|)}
&=&\cos{(\epp(t,\xi))}\left(\frac{t\sqrt{4\alpha^2 |\xi|^2-1}}{2\alpha}-t|\xi| \right) \nonumber \\[2mm]
&=&\cos{(\epp(t,\xi))} \cdot \frac{-t}{2\alpha \left(2\alpha |\xi|+\sqrt{4\alpha^2 |\xi|^2-1} \right)} \nonumber \\[2mm]
&=&  -\frac{t\cdot \cos{(\epp(t,\xi))}}{4\alpha^2 |\xi| \left(1+\sqrt{1-\frac{1}{4\alpha^2 |\xi|^2}}\right)}, \nonumber \\[2mm]
\cos{\left(\frac{t\sqrt{4\alpha^2|\xi|^2-1}}{2\alpha}\right)}-\cos{(t|\xi|)}&=& \frac{t\cdot \sin{(\eta(t,\xi))}}{4\alpha^2 |\xi| \left(1+\sqrt{1-\frac{1}{4\alpha^2 |\xi|^2}}\right)}. \nonumber
\end{eqnarray}
Here, $\epp(t,\xi)$ and $\eta(t,\xi)$ are some functions satisfying 
\[
\frac{t\sqrt{4\alpha^2 |\xi|^2-1}}{2\alpha}=t\sqrt{|\xi|^2-\frac{1}{4\alpha^2}}<\epp(t,\xi)\,,\,\eta(t,\xi)<t|\xi|.
\]
We also see that 
\begin{eqnarray}
\frac{1}{\frac{\sqrt{4\alpha^2 |\xi|^2-1}}{2\alpha}}-\frac{1}{|\xi|}
=\frac{2\alpha |\xi|-\sqrt{4\alpha^2 |\xi|^2-1}}{|\xi| \sqrt{4\alpha^2 |\xi|^2-1}}
&=&\frac{1}{|\xi| \sqrt{4\alpha^2 |\xi|^2-1}\, (2\alpha |\xi|+\sqrt{4\alpha^2 |\xi|^2-1})} \nonumber \\[2mm]
&=&\frac{1}{4\alpha^2 |\xi|^3 \sqrt{1-\frac{1}{4\alpha^2 |\xi|^2}} \left(1+\sqrt{1-\frac{1}{4\alpha^2 |\xi|^2}}\right)}. \nonumber
\end{eqnarray}
Together with (2.4), we rearrange 
\begin{eqnarray}
\hat{u}(t,\xi)&=&\hat{u}_{1}(\xi) \cdot E_{1}(t,\xi)+\hat{u}_{0}(\xi) \cdot E_{0}(t,\xi)+\hat{u}_{0}(\xi) \cdot \frac{1}{2\alpha} \cdot E_{1}(t,\xi) \nonumber \\[2mm]
&=&\hat{u}_{1}(\xi) \cdot e^{-\frac{t}{2\alpha}} \left[\frac{\sin{(t|\xi|)}}{|\xi|}-\frac{t \cdot \cos{(\epp(t,\xi))}}{4\alpha^2 |\xi|^2 \left(1+\sqrt{1-\frac{1}{4\alpha^2 |\xi|^2}}\right)}\right]  \nonumber \\[2mm] 
& &+\hat{u}_{1}(\xi) \cdot e^{-\frac{t}{2\alpha}} \sin{\left(\frac{t\sqrt{4\alpha^2 |\xi|^2-1}}{2\alpha} \right)} \cdot \frac{1}{4\alpha^2 |\xi|^3 \sqrt{1-\frac{1}{4\alpha^2 |\xi|^2}} \left(1+\sqrt{1-\frac{1}{4\alpha^2 |\xi|^2}}\right)} \nonumber \\[2mm]
& &+\hat{u}_{0}(\xi) \cdot e^{-\frac{t}{2\alpha}} \left[\cos{(t|\xi|)}+\frac{t \cdot \sin{(\eta(t,\xi))}}{4\alpha^2 |\xi| \left(1+\sqrt{1-\frac{1}{4\alpha^2 |\xi|^2}}\right)} \right] 
+\hat{u}_{0}(\xi) \cdot \frac{1}{2\alpha} \cdot E_{1}(t,\xi). \nonumber
\end{eqnarray}
Integrating over $\{\xi\in{\bf R}^{n}:|\xi|\ge1\}$, it follows that 
\[ \int_{|\xi| \ge 1} \left|\hat{u}(t,\xi)-\left(\hat{u}_{1}(\xi) \cdot e^{-\frac{t}{2|\xi|^2}} \frac{\sin{(t|\xi|)}}{|\xi|}+\hat{u}_{0}(\xi) \cdot e^{-\frac{t}{2|\xi|^2}} \cos{(t|\xi|)} \right) \right|^2 d\xi \]
\begin{eqnarray}
&\le& C\int_{|\xi| \ge 1} \left|\hat{u}_{1}(\xi) \cdot \left(e^{-\frac{t}{2\alpha}}-e^{-\frac{t}{2|\xi|^2}} \right) \cdot \frac{\sin{(t|\xi|)}}{|\xi|} \right|^2 d\xi \nonumber \\[2mm] 
& &+C\int_{|\xi| \ge 1} \left|\hat{u}_{0}(\xi) \cdot \left(e^{-\frac{t}{2\alpha}}-e^{-\frac{t}{2|\xi|^2}} \right) \cdot \cos{(t|\xi|)} \right|^2 d\xi \nonumber \\[2mm]
& &+C\int_{|\xi| \ge 1} \left|\hat{u}_{1}(\xi) \cdot e^{-\frac{t}{2\alpha}} \frac{t \cdot \cos{(\epp(t,\xi))}}{4\alpha^2 |\xi|^2 \left(1+\sqrt{1-\frac{1}{4\alpha^2 |\xi|^2}}\right)} \right|^2 d\xi \nonumber
\end{eqnarray}
\[
+C\int_{|\xi| \ge 1}\left|\hat{u}_{1}(\xi) \cdot e^{-\frac{t}{2\alpha}} \sin{\left(\frac{t\sqrt{4\alpha^2 |\xi|^2-1}}{2\alpha} \right)} \cdot \frac{1}{4\alpha^2 |\xi|^3 \sqrt{1-\frac{1}{4\alpha^2 |\xi|^2}} \left(1+\sqrt{1-\frac{1}{4\alpha^2 |\xi|^2}}\right)} \right|^2 d\xi
\]
\[
+C\int_{|\xi| \ge 1} \left|\hat{u}_{0}(\xi) \cdot e^{-\frac{t}{2\alpha}} \frac{t \cdot \sin{(\eta(t,\xi))}}{4\alpha^2 |\xi|\,\left(1+\sqrt{1-\frac{1}{4\alpha^2 |\xi|^2}}\right)} \right|^2 d\xi 
+C\int_{|\xi| \ge 1} \left|\hat{u}_{0}(\xi) \cdot \frac{1}{2\alpha} \cdot E_{1}(t,\xi) \right|^2 d\xi. \]
The first term and the second term on the right-hand side above can be estimated as follows:
\begin{eqnarray}
& &\int_{|\xi| \ge 1} \left|\hat{u}_{1}(\xi) \cdot \left(e^{-\frac{t}{2\alpha}}-e^{-\frac{t}{2|\xi|^2}} \right) \cdot \frac{\sin{(t|\xi|)}}{|\xi|} \right|^2 d\xi \nonumber \\[2mm]
&\le& \int_{|\xi| \ge 1} |\hat{u}_{1}(\xi)|^2 \cdot e^{-\frac{t}{\alpha}} \left|1-e^{-\frac{t}{2\alpha|\xi|^2}} \right|^2 \cdot \frac{1}{|\xi|^2} d\xi \nonumber \\[2mm]
&\le& \int_{|\xi| \ge 1} |\hat{u}_{1}(\xi)|^2 \cdot e^{-\frac{t}{\alpha}} \left(\frac{t}{2\alpha|\xi|^2} \right)^2 \cdot \frac{1}{|\xi|^2} d\xi \nonumber \\[2mm]
&\le& C(1+t)^2 \int_{|\xi| \ge 1}  e^{-\frac{t}{\alpha}} \cdot \frac{1}{\alpha^{l+2}\,|\xi|^6} \cdot \left(1+|\xi|^2 \right)^{l}|\hat{u}_{1}(\xi)|^2  d\xi \nonumber \\[2mm]
&\le& C(1+t)^2 \int_{|\xi| \ge 1} \frac{e^{-\frac{t}{\alpha}}}{\alpha^{l+5}} \cdot \left(1+|\xi|^2 \right)^{l}|\hat{u}_{1}(\xi)|^2  d\xi 
\le C \|u_{1}\|_{H^l}^2\, (1+t)^{-l-3}\quad(t\ge0), \nonumber
\end{eqnarray} 
\begin{eqnarray}
\int_{|\xi| \ge 1} \left|\hat{u}_{0}(\xi) \cdot \left(e^{-\frac{t}{2\alpha}}-e^{-\frac{t}{2|\xi|^2}} \right) \cdot \cos{(t|\xi|)} \right|^2 d\xi
&\le& \int_{|\xi| \ge 1} |\hat{u}_{0}(\xi)|^2 \cdot e^{-\frac{t}{\alpha}} \left(\frac{t}{2\alpha|\xi|^2} \right)^2 d\xi \nonumber \\[2mm]
&\le& C \|u_{0}\|_{H^{l+1}}^2\, (1+t)^{-l-3} 
\quad 
(t\ge0).
\nonumber
\end{eqnarray}
One can expect the remaining four terms decay faster than these two terms. 
Indeed, we confirm 
\begin{eqnarray}
\int_{|\xi| \ge 1} \left|\hat{u}_{1}(\xi) \cdot e^{-\frac{t}{2\alpha}} \frac{t \cdot \cos{(\epp(t,\xi))}}{4\alpha^2 |\xi|^2 \left(1+\sqrt{1-\frac{1}{4\alpha^2 |\xi|^2}}\right)} \right|^2 d\xi
&\le& C(1+t)^2 \int_{|\xi| \ge 1} |\hat{u}_{1}(\xi)|^2 \cdot e^{-\frac{t}{\alpha}} \cdot \frac{1}{\alpha^4 |\xi|^4} d\xi \nonumber \\[2mm]
&\le& C \|u_{1}\|_{H^l}^2\, (1+t)^{-l-4}
\quad 
(t\ge0),
\nonumber
\end{eqnarray}
\begin{eqnarray}
& &\int_{|\xi| \ge 1}\left|\hat{u}_{1}(\xi) \cdot e^{-\frac{t}{2\alpha}} \sin{\left(\frac{t\sqrt{4\alpha^2 |\xi|^2-1}}{2\alpha} \right)} \cdot \frac{1}{4\alpha^2 |\xi|^3 \sqrt{1-\frac{1}{4\alpha^2 |\xi|^2}} \left(1+\sqrt{1-\frac{1}{4\alpha^2 |\xi|^2}}\right)} \right|^2 d\xi \nonumber \\[2mm]
&\le& C \int_{|\xi| \ge 1} |\hat{u}_{1}(\xi)|^2 \cdot e^{-\frac{t}{\alpha}} \cdot \frac{1}{\alpha^4 |\xi|^6} d\xi
\ \le\ C \|u_{1}\|_{H^l}^2\, (1+t)^{-l-7} 
\quad 
(t\ge0),
\nonumber
\end{eqnarray}
\begin{eqnarray}
\int_{|\xi| \ge 1} \left|\hat{u}_{0}(\xi) \cdot e^{-\frac{t}{2\alpha}} \frac{t \cdot \sin{(\eta(t,\xi))}}{4\alpha^2 |\xi| \left(1+\sqrt{1-\frac{1}{4\alpha^2 |\xi|^2}}\right)} \right|^2 d\xi
&\le& C(1+t)^2 \int_{|\xi| \ge 1} |\hat{u}_{0}(\xi)|^2 \cdot e^{-\frac{t}{\alpha}} \cdot \frac{1}{\alpha^4 |\xi|^2} d\xi \nonumber \\[2mm]
&\le& C \|u_{0}\|_{H^{l+1}}^2\, (1+t)^{-l-4} 
\quad 
(t\ge0),
\nonumber
\end{eqnarray}
and 
\begin{eqnarray}
\int_{|\xi| \ge 1} \left|\hat{u}_{0}(\xi) \cdot \frac{1}{2\alpha} \cdot E_{1}(t,\xi) \right|^2 d\xi
&\le& C \int_{|\xi| \ge 1} |\hat{u}_{0}(\xi)|^2  \cdot e^{-\frac{t}{\alpha}} \cdot \frac{1}{\alpha^2|\xi|^2} d\xi \nonumber \\[2mm]
&\le& C \|u_{0}\|_{H^{l+1}}^2\, (1+t)^{-l-4} \quad (t\ge0).
\nonumber
\end{eqnarray}
Therefore, one obtains
\begin{eqnarray}
& &\int_{|\xi| \ge 1} \left|\hat{u}(t,\xi)-\left(\hat{u}_{1}(\xi) \cdot e^{-\frac{t}{2|\xi|^2}}  \frac{\sin{(t|\xi|)}}{|\xi|}+\hat{u}_{0}(\xi) \cdot e^{-\frac{t}{2|\xi|^2}}\,\cos{(t|\xi|)} \right) \right|^2 d\xi \nonumber \\[2mm]
&\le& C\left(||u_{1}||_{H^l}^2+||u_{0}||_{H^{l+1}}^2\right)\, (1+t)^{-l-3} \quad (t \ge 0), \nonumber 
\end{eqnarray}
and the proof is now complete. 
\hfill
$\Box$ 
\par
\vspace{0.1cm}

{\it Proof of Lemma 4.2.} 
Recall (2.2) to derive 
\begin{eqnarray}
\hat{u}(t,\xi)
&=&\frac{e^{\lambda_{1}t}-e^{\lambda_{2}t}}{\lambda_{1}-\lambda_{2}}\hat{u}_{1}(\xi)+\frac{\lambda_{1}e^{\lambda_{2}t}-\lambda_{2}e^{\lambda_{1}t}}{\lambda_{1}-\lambda_{2}}\hat{u}_{0}(\xi) \nonumber \\[2mm]
&=&\hat{u}_{1}(\xi) \cdot \frac{e^{\lambda_{1}t}-e^{\lambda_{2}t}}{\sqrt{1-4|\xi|^2 (1+|\xi|)^2}}
+\hat{u}_{1}(\xi) \cdot \frac{|\xi|^2 (e^{\lambda_{1}t}-e^{\lambda_{2}t})}{\sqrt{1-4|\xi|^2 (1+|\xi|)^2}} \nonumber \\[2mm]
& &+\hat{u}_{0}(\xi) \cdot \frac{e^{\lambda_{1}t}-e^{\lambda_{2}t}}{2\sqrt{1-4|\xi|^2 (1+|\xi|)^2}}
+\hat{u}_{0}(\xi) \cdot \frac{e^{\lambda_{1}t}+e^{\lambda_{2}t}}{2} \nonumber \\[2mm]
&=&\hat{u}_{1}(\xi) \cdot e^{\lambda_{1}t}-\hat{u}_{1}(\xi) \cdot e^{\lambda_{2}t}
+\hat{u}_{1}(\xi) \cdot f(\xi)\, (e^{\lambda_{1}t}-e^{\lambda_{2}t}) 
+\hat{u}_{1}(\xi) \cdot \frac{|\xi|^2 (e^{\lambda_{1}t}-e^{\lambda_{2}t})}{\sqrt{1-4|\xi|^2 (1+|\xi|)^2}} \nonumber \\[2mm]
& &+\hat{u}_{0}(\xi) \cdot e^{\lambda_{1}t}
+\hat{u}_{0}(\xi) \cdot \frac{f(\xi)}{2}\, (e^{\lambda_{1}t}-e^{\lambda_{2}t}), \nonumber
\end{eqnarray}
where
\[f(\xi):=\frac{1}{\sqrt{1-4|\xi|^2 (1+|\xi|)^2}}-1=\frac{4|\xi|^2 (1+|\xi|^2)^2}{\sqrt{1-4|\xi|^2 (1+|\xi|)^2}\left(1+\sqrt{1-4|\xi|^2 (1+|\xi|)^2}\right)}. \]
On the other hand, $e^{\lambda_{1}t}$ can be rewritten as 
\[ 
e^{\lambda_{1}t}=\exp\left(-t|\xi|^2 \cdot g(|\xi|^2) \right) 
\quad
\mbox{with}
\quad
g(\beta):=\frac{2(1+\beta)}{1+\sqrt{1-4\beta(1+\beta)^2}}.
 \]
The mean value theorem yields 
\begin{eqnarray}
e^{\lambda_{1}t}-e^{-t|\xi|^2}
&=&e^{-t|\xi|^2 \cdot g(|\xi|^2)}-e^{-t|\xi|^2} \nonumber \\[2mm]
&=&-t|\xi|^4 \, e^{-t|\xi|^2 \cdot g(\theta|\xi|^2)} 
\cdot \frac{dg}{d\beta}(\theta|\xi|^2) \nonumber
\end{eqnarray}
for some $0<\theta<1$. 
Since
\[
\sqrt{1-4\beta(1+\beta)^2}\ge\frac{1}{2}
\qquad
\mbox{for}
\quad
0 \le \beta \le \delta,
\]
there exist constants $c>0$ and $C>0$ such that 
\[
g(\beta)\ge c,
\qquad
\left|
\frac{dg}{d\beta}(\beta)
\right|
\le C,
\]
for $0\le\beta\le\delta$. 
for all $\xi\in{\bf R}^{n}$ satisfying $|\xi|\le\delta$. 
Therefore one can further calculate 
\begin{eqnarray}
\hat{u}(t,\xi)
&=&\hat{u}_{1}(\xi) \cdot e^{\lambda_{1}t}-\hat{u}_{1}(\xi) \cdot e^{\lambda_{2}t}
+\hat{u}_{1}(\xi) \cdot f(\xi)\, (e^{\lambda_{1}t}-e^{\lambda_{2}t}) 
+\hat{u}_{1}(\xi) \cdot \frac{|\xi|^2 (e^{\lambda_{1}t}-e^{\lambda_{2}t})}{\sqrt{1-4|\xi|^2 (1+|\xi|)^2}} \nonumber \\[2mm]
& &+\hat{u}_{0}(\xi) \cdot e^{\lambda_{1}t}
+\hat{u}_{0}(\xi) \cdot \frac{f(\xi)}{2}\, (e^{\lambda_{1}t}-e^{\lambda_{2}t}) \nonumber \\[2mm]
&=&P_{1} \cdot e^{-t|\xi|^2}+(A_{1}-iB_{1}) \, e^{-t|\xi|^2}-\hat{u}_{1}(\xi) \cdot t|\xi|^4\, e^{-t|\xi|^2 \cdot g(\theta|\xi|^2)} \cdot \frac{dg}{d\beta}(\theta|\xi|^2) \nonumber \\[2mm]
& &-\hat{u}_{1}(\xi) \cdot e^{\lambda_{2}t}
+\hat{u}_{1}(\xi) \cdot f(\xi) \, (e^{\lambda_{1}t}-e^{\lambda_{2}t}) 
+\hat{u}_{1}(\xi) \cdot \frac{|\xi|^2 (e^{\lambda_{1}t}-e^{\lambda_{2}t})}{\sqrt{1-4|\xi|^2 (1+|\xi|)^2}} \nonumber \\[2mm]
& &+P_{0} \cdot e^{-t|\xi|^2}+(A_{0}-iB_{0}) \, e^{-t|\xi|^2}-\hat{u}_{0}(\xi) \cdot t|\xi|^4 \, e^{-t|\xi|^2 \cdot g(\theta|\xi|^2)} \cdot \frac{dg}{d\beta}(\theta|\xi|^2) \nonumber \\[2mm]
& &+\hat{u}_{0}(\xi) \cdot \frac{f(\xi)}{2} \, (e^{\lambda_{1}t}-e^{\lambda_{2}t}). \nonumber
\end{eqnarray}
So one arrives at 
\begin{eqnarray}
& &\int_{|\xi| \le \delta} \left|\hat{u}(t,\xi)-(P_{1}+P_{0}) \, e^{-t|\xi|^2} \right|^2 d\xi \nonumber \\[2mm]
&\le&C\int_{|\xi| \le \delta} \left|(A_{1}-iB_{1}) \, e^{-t|\xi|^2} \right|^2 d\xi
+C\int_{|\xi| \le \delta} \left|\hat{u}_{1}(\xi) \cdot t|\xi|^4 \, e^{-t|\xi|^2 \cdot g(\theta|\xi|^2)} \cdot \frac{dg}{d\beta}(\theta|\xi|^2) \right|^2 d\xi \nonumber \\[2mm]
& &+C\int_{|\xi| \le \delta} \left|\hat{u}_{1}(\xi) \cdot e^{\lambda_{2}t} \right|^2 d\xi
+C\int_{|\xi| \le \delta} \left|\hat{u}_{1}(\xi) \cdot f(\xi) \, (e^{\lambda_{1}t}-e^{\lambda_{2}t}) \right|^2 d\xi \nonumber \\[2mm]
& &+C\int_{|\xi| \le \delta} \left|\hat{u}_{1}(\xi) \cdot \frac{|\xi|^2 (e^{\lambda_{1}t}-e^{\lambda_{2}t})}{\sqrt{1-4|\xi|^2 (1+|\xi|)^2}} \right|^2 d\xi
+C\int_{|\xi| \le \delta} \left|(A_{0}-iB_{0}) \, e^{-t|\xi|^2} \right|^2 d\xi \nonumber \\[2mm]
& &+C\int_{|\xi| \le \delta} \left|\hat{u}_{0}(\xi) \cdot t|\xi|^4 \, e^{-t|\xi|^2 \cdot g(\theta|\xi|^2)} \cdot \frac{dg}{d\beta}(\theta|\xi|^2) \right|^2 d\xi \nonumber \\[2mm]
& &+C\int_{|\xi| \le \delta} \left|\hat{u}_{0}(\xi) \cdot \frac{f(\xi)}{2} \, (e^{\lambda_{1}t}-e^{\lambda_{2}t})\right|^2d\xi. \nonumber
\end{eqnarray}

Let us evaluate each term in the right hand side of the inequality just above:
\begin{eqnarray}
\int_{|\xi| \le \delta} \left|(A_{1}-iB_{1}) \, e^{-t|\xi|^2} \right|^2 d\xi
&\le& C \|u_{1}\|^2_{1,1} \int_{|\xi| \le \delta} |\xi|^2 \, e^{-2t|\xi|^2} d\xi \nonumber \\[2mm]
&\le& C \|u_{1}\|^2_{1,1}\, (1+t)^{-\frac{n+2}{2}} \quad (t \ge 0); \nonumber
\end{eqnarray}
\begin{eqnarray}
\int_{|\xi| \le \delta} \left|\hat{u}_{1} \cdot t|\xi|^4 \, e^{-t|\xi|^2 \cdot g(\theta|\xi|^2)} \cdot \frac{dg}{d\beta}(\theta|\xi|^2) \right|^2 d\xi
&\le&C\|u_{1}\|^2_{1}\, (1+t)^2 \int_{|\xi| \le \delta} |\xi|^8 \, e^{-ct|\xi|^2} d\xi \nonumber \\[2mm]
&\le&C\|u_{1}\|^2_{1}\, (1+t)^{-\frac{n+4}{2}} \quad (t \ge 0); \nonumber
\end{eqnarray}
\begin{eqnarray}
\int_{|\xi| \le \delta} \left|\hat{u}_{1}(\xi) \cdot e^{\lambda_{2}t} \right|^2 d\xi
&\le&\int_{|\xi| \le \delta} |\hat{u}_{1}|^2 e^{-ct} d\xi
\le C\|u_{1}\|^2_{2} \, e^{-ct} \quad (t \ge 0); \nonumber
\end{eqnarray}
\begin{eqnarray}
\int_{|\xi| \le \delta} \left|\hat{u}_{1}(\xi) \cdot f(\xi) \, (e^{\lambda_{1}t}-e^{\lambda_{2}t}) \right|^2 d\xi 
&\le& C\|u_{1}\|^2_{1} \int_{|\xi| \le \delta} |\xi|^4 \, e^{2\lambda_{1}t} d\xi \nonumber \\[2mm]
&\le& C\|u_{1}\|^2_{1} \int_{|\xi| \le \delta} |\xi|^4 \, e^{-ct|\xi|^2} d\xi \nonumber \\[2mm]
&\le& C\|u_{1}\|^2_{1}\,(1+t)^{-\frac{n+4}{2}} \quad (t \ge 0); \nonumber
\end{eqnarray}
\begin{eqnarray}
\int_{|\xi| \le \delta} \left|\hat{u}_{1}(\xi) \cdot \frac{|\xi|^2 (e^{\lambda_{1}t}-e^{\lambda_{2}t})}{\sqrt{1-4|\xi|^2 (1+|\xi|)^2}} \right|^2 d\xi
&\le&C\|u_{1}\|^2_{1} \int_{|\xi| \le \delta} |\xi|^4 \, e^{2\lambda_{1}t} d\xi \nonumber \\[2mm]
&\le&C\|u_{1}\|^2_{1}\, (1+t)^{-\frac{n+4}{2}} \quad (t \ge 0); \nonumber
\end{eqnarray}
\begin{eqnarray}
\int_{|\xi| \le \delta} \left|(A_{0}-iB_{0}) \, e^{-t|\xi|^2} \right|^2 d\xi
&\le&C\|u_{0}\|^2_{1,1} \int_{|\xi| \le \delta} |\xi|^2 \, e^{-2t|\xi|^2} d\xi \nonumber \\[2mm]
&\le& C\|u_{1}\|^2_{1,1} (1+t)^{-\frac{n+2}{2}} \quad (t \ge 0); \nonumber
\end{eqnarray}
\begin{eqnarray}
\int_{|\xi| \le \delta} \left|\hat{u}_{0}(\xi) \cdot t|\xi|^4 \, e^{-t|\xi|^2 \cdot g(\theta|\xi|^2)} \cdot \frac{dg}{d\beta}(\theta|\xi|^2) \right|^2 d\xi
&\le&C\|u_{0}\|^2_{1}\, (1+t)^2 \int_{|\xi| \le \delta} |\xi|^8 \, e^{-ct|\xi|^2} d\xi \nonumber \\[2mm]
&\le&C\|u_{1}\|^2_{1}\, (1+t)^{-\frac{n+4}{2}} \quad (t \ge 0); \nonumber
\end{eqnarray}
\begin{eqnarray}
\int_{|\xi| \le \delta} \left|\hat{u}_{0}(\xi) \cdot \frac{f(\xi)}{2} \, (e^{\lambda_{1}t}-e^{\lambda_{2}t})\right|^2d\xi
&\le&C\|u_{0}\|^2_{1} \int_{|\xi| \le \delta} |\xi|^4 \, e^{2\lambda_{1}t} d\xi \nonumber \\[2mm]
&\le&C\|u_{1}\|^2_{1}\, (1+t)^{-\frac{n+4}{2}} \quad (t \ge 0). \nonumber
\end{eqnarray}
Therefore, one obtains 
\begin{eqnarray}
\int_{|\xi| \le \delta} \left|\hat{u}(t,\xi)-(P_{1}+P_{0}) \, e^{-t|\xi|^2} \right|^2 d\xi 
&\le& C\left(\|u_{1}\|_{1,1}^2+\|u_{0}\|_{1,1}^2\right) \, (1+t)^{-\frac{n+2}{2}} \quad (t \ge 0), \nonumber
\end{eqnarray}
which is the desired estimate. 
\hfill
$\Box$
\par
We can also obtain corresponding estimates in the middle frequency region with the aid of the following lemmas. 
\begin{lem}
It holds 
\[\sup_{x > 0}
\left|
{\frac{\sin{tx}}{x}}
\right|
 = t \]
 for $t \ge 0$. 
\end{lem}
\begin{lem}
There exists a constant $C>0$ such that 
\[\frac{\sinh{tx}}{x} \le Cte^{tx} \]
for $x>0$ and $t \ge 0$.
\end{lem}
Proofs of these lemmas above are elementary, and so omitted.\\ 

{\it Proof of Lemma 4.3.}\,By using $\zeta$ defined in (2.4), one first sees 
\begin{eqnarray}
\int_{\delta<|\xi|<1} |\hat{u}(t,\xi)|^2 d\xi
&\le&C\int_{\delta<|\xi| \le \zeta} \left|\hat{u}_{0}(\xi) \cdot E_{0}(t,\xi) \right|^2 d\xi
+C\int_{\zeta<|\xi|<1} \left|\hat{u}_{0}(\xi) \cdot E_{0}(t,\xi) \right|^2 d\xi \nonumber \\[2mm]
& &+C\int_{\delta<|\xi| \le \zeta} \left|\left\{\hat{u}_{1}(\xi)+\frac{1}{2(1+|\xi|^2)} \cdot \hat{u}_{0}(\xi) \right\} E_{1}(t,\xi) \right|^2 d\xi \nonumber \\[2mm]
& &+C\int_{\zeta<|\xi|<1} \left|\left\{\hat{u}_{1}(\xi)+\frac{1}{2(1+|\xi|^2)} \cdot \hat{u}_{0}(\xi) \right\} E_{1}(t,\xi) \right|^2 d\xi. \nonumber
\end{eqnarray}
We evaluate of each section above by using Lemmas 4.6 and 4.7. Indeed, there exists a constant $\eta > 0$ such that  
\begin{eqnarray}
\int_{\delta<|\xi| \le \zeta} \left|\hat{u}_{0}(\xi) \cdot E_{0}(t,\xi) \right|^2 d\xi
&\le& \int_{\delta<|\xi| \le \zeta} |\hat{u}_{0}(\xi)|^2 \cdot e^{-\frac{t}{(1+|\xi|^2)}} \cdot e^{\frac{t\sqrt{1-4|\xi|^2(1+|\xi|^2)^2}}{(1+|\xi|^2)}} d\xi \nonumber \\[2mm]
&\le& \int_{\delta<|\xi| \le \zeta} |\hat{u}_{0}(\xi)|^2 \cdot e^{-\eta t} d\xi
 \le \|u_{0}\|_{2}^2 \, e^{-\eta t}, \nonumber
\end{eqnarray}
\begin{eqnarray}
& &\int_{\delta<|\xi| \le \zeta} \left|\left\{\hat{u}_{1}(\xi)+\frac{1}{2(1+|\xi|^2)} \cdot \hat{u}_{0}(\xi) \right\} E_{1}(t,\xi) \right|^2 d\xi \nonumber \\[2mm]
&\le& C \int_{\delta<|\xi| \le \zeta} \left(|\hat{u}_{1}(\xi)|^2+|\hat{u}_{0}(\xi)|^2 \right) \left|E_{1}(t,\xi)\right|^2 d\xi \nonumber \\[2mm]
&\le& Ct^2 \int_{\delta<|\xi| \le \zeta} \left(|\hat{u}_{1}(\xi)|^2+|\hat{u}_{0}(\xi)|^2 \right) e^{-\frac{t}{(1+|\xi|^2)}} \cdot e^{\frac{t\sqrt{1-4|\xi|^2(1+|\xi|^2)^2}}{(1+|\xi|^2)}} d\xi \nonumber \\[2mm]
&\le& Ct^2 \left(\|u_{0}\|_{2}^2+\|u_{1}\|_{2}^2 \right)\, e^{-\eta t}
 \le C\left(\|u_{0}\|_{2}^2+\|u_{1}\|_{2}^2 \right) \, e^{-\eta' t},
 \quad
 0<\eta'<\eta,
 \nonumber
\end{eqnarray}
and
\begin{eqnarray}
\int_{\zeta<|\xi|<1} \left|\hat{u}_{0}(\xi) \cdot E_{0}(t,\xi) \right|^2 d\xi
&\le&\int_{\zeta<|\xi|<1} |\hat{u}_{0}(\xi)|^2 \cdot e^{-\frac{t}{(1+|\xi|^2)}} d\xi
\ \le\ \|u_{0}\|_{2}^2 \, e^{-\eta t}, \nonumber
\end{eqnarray}
\begin{eqnarray}
& &\int_{\zeta<|\xi|<1} \left|\left\{\hat{u}_{1}(\xi)+\frac{1}{2(1+|\xi|^2)} \cdot \hat{u}_{0}(\xi) \right\} E_{1}(t,\xi) \right|^2 d\xi \nonumber \\[2mm]
&\le& Ct^2\int_{\zeta<|\xi|<1} \left(|\hat{u}_{1}(\xi)|^2+|\hat{u}_{0}(\xi)|^2 \right) 
e^{-\frac{t}{(1+|\xi|^2)}} d\xi
\le C\left(\|u_{0}\|_{2}^2+\|u_{1}\|_{2}^2 \right) \, e^{-\eta' t}. \nonumber
\end{eqnarray}
Therefore, one obtains the desired inequality.
\hfill
$\Box$ 
\par
Now, we check that 
\begin{eqnarray}
& &\int_{|\xi| \le 1} \left|\hat{u}_{1}(\xi) \cdot e^{-\frac{t}{2|\xi|^2}} \frac{\sin{(t|\xi|)}}{|\xi|}+\hat{u}_{0}(\xi) \cdot e^{-\frac{t}{2|\xi|^2}} \cos{(t|\xi|)} \right|^2 d\xi \nonumber \\[2mm]
&\le& C\int_{|\xi| \le 1} \left|\hat{u}_{1}(\xi) \cdot e^{-\frac{t}{2|\xi|^2}} \frac{\sin{(t|\xi|)}}{|\xi|} \right|^2 d\xi
+C\int_{|\xi| \le 1} \left|\hat{u}_{0}(\xi) \cdot e^{-\frac{t}{2|\xi|^2}} \cos{(t|\xi|)} \right|^2 d\xi \nonumber \\[2mm]
&\le& Ct^2 \int_{|\xi| \le 1} |\hat{u}_{1}(\xi)|^2 \cdot e^{-\frac{t}{|\xi|^2}} d\xi
+C\int_{|\xi| \le 1} |\hat{u}_{0}(\xi)|^2 \cdot e^{-\frac{t}{|\xi|^2}} d\xi \nonumber \\[2mm]
&\le& Ct^2 e^{-t} \int_{|\xi| \le 1} |\hat{u}_{1}(\xi)|^2 d\xi
+C e^{-t} \int_{|\xi| \le 1} |\hat{u}_{0}(\xi)|^2 d\xi \nonumber \\[2mm]
&\le& C(t^2\|u_1\|_2^2+\|u_0\|_2^2) \, e^{-t}
\le C(\|u_1\|_2^2+\|u_0\|_2^2) \, e^{-\eta t}
\end{eqnarray}
for $t \ge 0$. Here, $0<\eta<1$ is a constant independent of $t$ and $u_0$, $u_1$.
From (4.1) and Lemmas 4.1-4.3, it follows that  
\[\int_{{\bf R}^{n}} \left|\hat{u}(t,\xi)-(P_1+P_0)\, e^{-t|\xi|^2}-\left(\hat{u}_{1}(\xi) \cdot e^{-\frac{t}{2|\xi|^2}} \frac{\sin{(t|\xi|)}}{|\xi|}+\hat{u}_{0}(\xi) \cdot e^{-\frac{t}{2|\xi|^2}} \cos{(t|\xi|)} \right) \right|^2 d\xi\]
\[= \int_{|\xi| \le \delta} \left|\hat{u}(t,\xi)-(P_1+P_0)\, e^{-t|\xi|^2}-\left(\hat{u}_{1}(\xi) \cdot e^{-\frac{t}{2|\xi|^2}} \frac{\sin{(t|\xi|)}}{|\xi|}+\hat{u}_{0}(\xi) \cdot e^{-\frac{t}{2|\xi|^2}} \cos{(t|\xi|)} \right) \right|^2 d\xi\]
\[+ \int_{\delta \le |\xi| \le 1} \left|\hat{u}(t,\xi)-(P_1+P_0)\, e^{-t|\xi|^2}-\left(\hat{u}_{1}(\xi) \cdot e^{-\frac{t}{2|\xi|^2}} \frac{\sin{(t|\xi|)}}{|\xi|}+\hat{u}_{0}(\xi) \cdot e^{-\frac{t}{2|\xi|^2}} \cos{(t|\xi|)} \right) \right|^2 d\xi\]
\[+ \int_{1 \le |\xi|} \left|\hat{u}(t,\xi)-(P_1+P_0)\, e^{-t|\xi|^2}-\left(\hat{u}_{1}(\xi) \cdot e^{-\frac{t}{2|\xi|^2}} \frac{\sin{(t|\xi|)}}{|\xi|}+\hat{u}_{0}(\xi) \cdot e^{-\frac{t}{2|\xi|^2}} \cos{(t|\xi|)} \right) \right|^2 d\xi\] 
\begin{eqnarray}
&\le& C\left(\|u_{1}\|_{1,1}^2+\|u_{0}\|_{1,1}^2\right) \, t^{-\frac{n+2}{2}}+C\left(\|u_{1}\|_2^2+\|u_{0}\|_2^2 \right) \, e^{-\eta t} \nonumber \\[2mm]
& &+ C \left(\|u_{0}\|_2^2+\|u_{1}\|_2^2 \right) \, e^{-\eta t}+C|P_1+P_0|^2 \, e^{-\eta t}+C\left(\|u_{1}\|_2^2+\|u_{0}\|_2^2 \right) \, e^{-\eta t} \nonumber \\[2mm]
& &+ C\left(||u_{1}||_{H^l}^2+||u_{0}||_{H^{l+1}}^2 \right) \, t^{-l-3}+C|P_1+P_0|^2 \, e^{-\eta t} \nonumber \\[2mm]
&\le& C\left(\|u_{1}\|_{1,1}^2+\|u_{0}\|_{1,1}^2\right) \, t^{-\frac{n+2}{2}}
+C\left(||u_{1}||^2_{H^l}+||u_{0}||^2_{H^{l+1}}\right) \, t^{-l-3} \nonumber \\[2mm]
& &+ C\left(|P_1+P_0|^2+\|u_{1}\|_2^2+\|u_{0}\|_2^2 \right) \, e^{-\eta t} \quad (t \ge 1),
\end{eqnarray}
where we have used the fact that
\[
\int_{|\xi| \ge 1}e^{-2t|\xi|^2} d\xi
\le e^{-t} \int_{|\xi| \ge 1}e^{-|\xi|^2} d\xi
\le Ce^{-t} \quad (t \ge 1).
\]
Here, $C>0$ is a constant independent of $t$.\\

By using (4.2), one can prove Theorems 3.1, 3.2 and 3.3.

{\it Proof of Theorem 3.1.}\,By relying on (4.2), one can proceed the following estimates: ($t \ge 1$)
\begin{eqnarray}
& &\left\|\hat{u}(t,\xi)-\left(\hat{u}_{1}(\xi) \cdot e^{-\frac{t}{2|\xi|^2}} \frac{\sin{(t|\xi|)}}{|\xi|}+\hat{u}_{0}(\xi) \cdot e^{-\frac{t}{2|\xi|^2}} \cos{(t|\xi|)} \right) \right\|_2
\ \le\ \left\|(P_0+P_1) \, e^{-t|\xi|^2} \right\|_2 \nonumber \\[2mm]
& &+\left\|\hat{u}(t,\xi)-(P_1+P_0) \, e^{-t|\xi|^2}-\left(\hat{u}_{1}(\xi) \cdot e^{-\frac{t}{2|\xi|^2}}\, \frac{\sin{(t|\xi|)}}{|\xi|}+\hat{u}_{0}(\xi) \cdot e^{-\frac{t}{2|\xi|^2}}\,\cos{(t|\xi|)} \right) \right\|_2 \nonumber \\[2mm]
&\le& C|P_{0}+P_{1}| \, t^{-\frac{n}{4}}+C\left(\|u_{1}\|_{1,1}+\|u_{0}\|_{1,1}\right) \, t^{-\frac{n+2}{4}}\nonumber \\[2mm]
& &+ C\left(||u_{1}||_{H^l}+||u_{0}||_{H^{l+1}}\right) \, t^{-\frac{l+3}{2}}
+ C\left(|P_1+P_0|+\|u_{1}\|_2+\|u_{0}\|_2 \right) \, e^{-\eta t} \nonumber \\[2mm]
&\le& \left\{ \begin{array}{ll}
C I_0 \, t^{-\frac{l+3}{2}}, 
& \quad 
l \leq n/2-3, \\[2mm]
C I_0 \, t^{-\frac{n}{4}},
& \quad 
n/2-3 <l < n/2-1.
\end{array}
\right.\nonumber
\end{eqnarray}
Taking into account of $l\ge2$, the proof is complete.
\hfill
$\Box$
\par 
{\it Proof of Theorem 3.2.}\, One also uses (4.2) to prove Theorem 3.2 as follows: ($t \ge 1$)
\begin{eqnarray}
& &\left\|\hat{u}(t,\xi)-(P_1+P_0) \, e^{-t|\xi|^2} \right\|_2
\ \le\ \left\|\hat{u}_{1}(\xi) \cdot e^{-\frac{t}{2|\xi|^2}} \frac{\sin{(t|\xi|)}}{|\xi|}+\hat{u}_{0}(\xi) \cdot e^{-\frac{t}{2|\xi|^2}} \cos{(t|\xi|)} \right\|_2 \nonumber \\[2mm]
& &+\left\|\hat{u}(t,\xi)-(P_1+P_0) \, e^{-t|\xi|^2}-\left(\hat{u}_{1}(\xi) \cdot e^{-\frac{t}{2|\xi|^2}} \frac{\sin{(t|\xi|)}}{|\xi|}+\hat{u}_{0}(\xi) \cdot e^{-\frac{t}{2|\xi|^2}} \cos{(t|\xi|)} \right) \right\|_2 \nonumber \\[2mm]
&\le& \left\|\hat{u}_{1}(\xi) \cdot e^{-\frac{t}{2|\xi|^2}} \frac{\sin{(t|\xi|)}}{|\xi|}+\hat{u}_{0}(\xi) \cdot e^{-\frac{t}{2|\xi|^2}} \cos{(t|\xi|)} \right\|_{|\xi| \le 1} \nonumber \\[2mm]
& &+\left\|\hat{u}_{1}(\xi) \cdot e^{-\frac{t}{2|\xi|^2}} \frac{\sin{(t|\xi|)}}{|\xi|}+\hat{u}_{0}(\xi) \cdot e^{-\frac{t}{2|\xi|^2}} \cos{(t|\xi|)} \right\|_{|\xi| \ge 1} \nonumber \\[2mm]
& &+\left\|\hat{u}(t,\xi)-(P_1+P_0) \, e^{-t|\xi|^2}-\left(\hat{u}_{1}(\xi) \cdot e^{-\frac{t}{2|\xi|^2}} \frac{\sin{(t|\xi|)}}{|\xi|}+\hat{u}_{0}(\xi) \cdot e^{-\frac{t}{2|\xi|^2}} \cos{(t|\xi|)} \right) \right\|_2 \nonumber \\[2mm]
&\le& C\left(||u_{1}||_{H^l}+||u_{0}||_{H^{l+1}} \right) \, t^{-\frac{l+1}{2}}
+C\left(\|u_{1}\|_{1,1}+\|u_{0}\|_{1,1} \right) \, t^{-\frac{n+2}{4}} \nonumber \\[2mm]
& &+C\left(||u_{1}||_{H^l}+||u_{0}||_{H^{l+1}}\right) \, t^{-\frac{l+3}{2}}
+C\left(|P_1+P_0|+\|u_{1}\|_2+\|u_{0}\|_2 \right) \, e^{-\eta t} \nonumber \\[2mm]
&\le& \left\{ \begin{array}{ll}
C I_0 \, t^{-\frac{l+1}{2}}, 
& \quad n/2-1<l \leq n/2, \\[2mm]
C I_0 \, t^{-\frac{n+2}{4}}, 
&\quad n/2<l.
\end{array}
\right.\nonumber
\end{eqnarray}
This implies the desired estimate.
\hfill
$\Box$
 
\vspace{0.2cm}
Finally, the proof of Theorem 3.3 is a direct consequence of (4.2). \\

\section{Decay estimate of solutions}
Based on Theorems 3.1 and 3.2, one can obtain optimal decay estimates of the $L^{2}$-norm of solutions to problem (1.1). Nowadays, under the results in Theorems 3.1-3.3, the following computations are standard by relying on the Plancherel theorem (cf. \cite{Ike-1, IS}). 

It follows from Theorems 3.1 and 3.2 that
\begin{eqnarray}
& &\|\hat{u}(t,\xi)\|_2
\ \le\ \left\|\hat{u}_{1}(\xi) \cdot e^{-\frac{t}{2|\xi|^2}} \frac{\sin{(t|\xi|)}}{|\xi|}+\hat{u}_{0}(\xi) \cdot e^{-\frac{t}{2|\xi|^2}} \cos{(t|\xi|)} \right\|_2 \nonumber \\[2mm]
& &+\left\|\hat{u}(t,\xi)-\left(\hat{u}_{1}(\xi) \cdot e^{-\frac{t}{2|\xi|^2}} \frac{\sin{(t|\xi|)}}{|\xi|}+\hat{u}_{0}(\xi) \cdot e^{-\frac{t}{2|\xi|^2}} \cos{(t|\xi|)} \right) \right\|_2 \nonumber \\[2mm]
&\le& \left\{ \begin{array}{ll}
C\left(||u_{1}||_{H^l}+||u_{0}||_{H^{l+1}} \right) \, t^{-\frac{l+1}{2}}
+CI_{0}\,t^{-\frac{l+3}{2}}, 
& \qquad 
2\le l\leq n/2-3 
\,\,\,
\mbox{with}
\,\,\,
n\ge10, \\[2mm]
C\left(||u_{1}||_{H^l}+||u_{0}||_{H^{l+1}} \right) \, t^{-\frac{l+1}{2}}
+CI_{0}\,t^{-\frac{n}{4}}, 
& \qquad 
2\le l <n/2 -1 
\,\,\,
\mbox{with}
\,\,\,
n=7, 8, 9 \\[2mm]
& \hspace{0.2cm}\mbox{or}\,\,\, 
n/2-3 <l < n/2-1
\,\,\,\mbox{with}\,\,\,n\ge10, \\[0.2cm]
\end{array}
\right. \nonumber 
\end{eqnarray}
for $t \ge 1$. Similarly to the above computation one has
\begin{eqnarray}
& & \|\hat{u}(t,\xi)\|_2
\le \left\|(P_0+P_1) \, e^{-t|\xi|^2}\right\|_2+\left\|\hat{u}(t,\xi)-(P_1+P_0) \, e^{-t|\xi|^2} \right\|_2 \nonumber \\[2mm]
&\le&  \left\{ \begin{array}{ll}
C|P_{0}+P_{1}| \, t^{-\frac{n}{4}}
+CI_{0}\,t^{-\frac{l+1}{2}}, 
&  \quad  
\ell=2 \,\,\,\mbox{with}\,\,\,n=4
\quad\mbox{or}\quad
2\le l\le 5/2 \,\,\,\mbox{with}\,\,\,n=5\\[2mm]
& \hspace{2cm}\mbox{or}\,\,\, 
n/2-1 <l \leq n/2
\,\,\,\mbox{with}\,\,\,n\ge6, \\[0.2cm]
C|P_{0}+P_{1}| \, t^{-\frac{n}{4}}
+CI_{0}\,t^{-\frac{n+2}{4}}, 
& \quad 
n/2<l\,\,\,\mbox{with}\,\,\,n\ge4\\[2mm]
& \hspace{2cm}\mbox{or}\,\,\, 
2 \leq l
\,\,\,\mbox{with}\,\,\,n = 1,2,3,\\[0.2cm]
\end{array}
\right. \nonumber 
\end{eqnarray}
for $t \ge 1$.
Furthermore, Theorem 3.3 leads us to the estimate in the case of $l=n/2-1$ with $n\ge6$:
\begin{eqnarray}
& &\|\hat{u}(t,\xi)\|_2 \nonumber \\[2mm] 
&\le& \left\|(P_0+P_1) \, e^{-t|\xi|^2}\right\|_2
+\left\|\hat{u}_{1}(\xi) \, e^{-\frac{t}{2|\xi|^2}} \frac{\sin{(t|\xi|)}}{|\xi|}+\hat{u}_{0}(\xi) \cdot e^{-\frac{t}{2|\xi|^2}} \cos{(t|\xi|)} \right\|_2 \nonumber \\[2mm]
& &+\left\|\hat{u}(t,\xi)-\left\{(P_1+P_0) \, e^{-t|\xi|^2}
+\hat{u}_{1}(\xi) \cdot e^{-\frac{t}{2|\xi|^2}} \frac{\sin{(t|\xi|)}}{|\xi|}+\hat{u}_{0}(\xi) \cdot e^{-\frac{t}{2|\xi|^2}}\cos{(t|\xi|)} \right\} \right\|_2 \nonumber \\[2mm]
&\le& C|P_{0}+P_{1}| \, t^{-\frac{n}{4}}
+C\left(||u_{1}||_{H^l}+||u_{0}||_{H^{l+1}} \right) \, t^{-\frac{l+1}{2}}+CI_{0}\,t^{-\frac{n+2}{4}} . \nonumber \\[2mm]
&\le& CI_{0}\,t^{-\frac{n}{4}} \nonumber
\end{eqnarray}
for $t \ge 1$. Thus, using the Plancherel theorem one has the following optimal decay results under the condition $l \geq 2$ to guarantee the (unique) existence of solutions to (1.1). 
\begin{theo}
\label{thm:4}
Let $n \ge 7$ and  $l$ satisfy $2 \leq l < n/2-1$.\\ 
If $(u_{0},u_{1}) \in (H^{l+1}({\bf R}^{n}) \cap L^{1,1}({\bf R}^{n})) \times (H^{l}({\bf R}^{n}) \cap L^{1,1}({\bf R}^{n}))$, 
then the solution $u \in X_{2}$ to {\rm (1.1)} satisfies 
\[\Vert u(t,\cdot)\Vert_2 \leq CI_{0}\,t^{-\frac{l+1}{2}} \quad (t \ge 1), \]
where $C>0$ is a constant independent of $t$ and $u_0,\,u_1$.
\end{theo}
\begin{theo}
\label{thm:5}
Let $n\ge1$ and $l$ satisfy 
\begin{eqnarray*}
\left\{
\begin{array}{lll}
2\le l, & \quad n=1,2,3,4,5, \\[2mm]
n/2-1 \leq l, &\quad n\ge6.
\end{array}
\right.
\end{eqnarray*}
If $(u_{0},u_{1}) \in (H^{l+1}({\bf R}^{n}) \cap L^{1,1}({\bf R}^{n})) \times (H^{l}({\bf R}^{n}) \cap L^{1,1}({\bf R}^{n}))$, 
then the solution $u \in X_{2}$ to {\rm (1.1)} satisfies 
\[\Vert u(t,\cdot)\Vert_2 \leq CI_{0}\,t^{-\frac{n}{4}} \quad (t \ge 1), \]
where $C>0$ is a constant independent of $t$ and $u_0,\,u_1$.
\end{theo}
\begin{rem}{\rm The number $l^{*} := n/2-1$ found in Theorems 5.1 and 5.2 is a kind of threshold that divides the property of the solution $u(t,x)$ to problem (1.1) into two parts: one is wave-like, and the other is parabolic-like. In this sense, the number $l^{*}$ seems to be quite meaningful.}
\end{rem}

For example, as a direct consequence of Theorems 3.1, 3.2, 3.3, 5.1 and 5.2, in the case of $n = 10$ one can obtain the following results.\,\, $l^{*} = 4$, and\\
{\rm (I)}\,$2 \leq l < 4$ $\Rightarrow$
\[\hat{u}(t,\xi) \sim \hat{u}_{1}(\xi) \, e^{-\frac{t}{2|\xi|^2}} \frac{\sin{(t|\xi|)}}{|\xi|}+\hat{u}_{0}(\xi) \cdot e^{-\frac{t}{2|\xi|^2}} \cos{(t|\xi|)}\quad(t \to \infty),\]  
{\rm (II)}\,$4 < l$ $\Rightarrow$
\[\hat{u}(t,\xi) \sim (P_0+P_1) \, e^{-t|\xi|^2}\,\,(t \to \infty),\]
and\\
{\rm (III)}\,$l = 4$ $\Rightarrow$
\[\hat{u}(t,\xi) \sim \hat{u}_{1}(\xi) \, e^{-\frac{t}{2|\xi|^2}} \frac{\sin{(t|\xi|)}}{|\xi|}+\hat{u}_{0}(\xi) \cdot e^{-\frac{t}{2|\xi|^2}} \cos{(t|\xi|)} + (P_0+P_1) \, e^{-t|\xi|^2} \,\,(t \to \infty).\]  
Furthermore, one can get the following optimal decay estimates of solutions in $L^{2}$-sense.\\ 
{\rm (IV)}\,$2 \leq l < 4$ $\Rightarrow$
\[\Vert u(t,\cdot)\Vert \leq CI_{0}t^{-\frac{l+1}{2}}\,\,(t \geq 1),\]
and\\
{\rm (V)}\,$4 \leq l$ $\Rightarrow$
\[\Vert u(t,\cdot)\Vert \leq CI_{0}t^{-\frac{5}{2}}\,\,(t \geq 1).\]

\par
\vspace{0.5cm}
\noindent{\em Acknowledgement.}
\smallskip
The work of the second author (R. IKEHATA) was supported in part by Grant-in-Aid for Scientific Research (C)15K04958  of JSPS. 


\end{document}